\documentclass[12pt,letterpaper]{amsart}

\usepackage{amssymb,latexsym,amsmath,amsthm,graphicx}
\numberwithin{equation}{section}

\newtheorem*{remark}{Remark}

\usepackage{xcolor}

\newcommand{\ol}{\overline}

\def \re{\text{Re}}

\input epsf
\def\ds{\displaystyle}
\def\M{\text{\rm mod}}

\begin{document}

\title{Remarks on the disjoint minimal graphs conjecture} 

\author[Allen Weitsman]{Allen Weitsman}

%\address{Department of Mathematics\\ Purdue University\\
%West Lafayette, IN 47907-1395}
\address{Email: weitsman@purdue.edu}

%\subjclass[2000]{Primary: , Secondary: }
%\date{2011}

\begin{abstract}
It is known that there can be at most 3 disjoint minimal graphs in $R^3$ with 0 boundary
values, and conjectured that there can be at most 2.  In this paper we explore 
properties that the graphs would have if there were 3.


\end{abstract}

\maketitle

%\tableofcontents
\section{Introduction} Let $\ D\ $ be an unbounded domain in $\Bbb R^2$.  
 We are interested in
solutions to the minimal surface equation  for the boundary value problem
\begin{equation}
\label{mineq1}
\begin{aligned}
Lu=\text{div}\left({\nabla u\over\sqrt{1+|\nabla u|^2}}\right)=0\\
u>0\quad\text{ in }\ \ D,\qquad  u=0\quad \text{ on }\  \partial D.
\end{aligned}
\end{equation} 
Henceforth we shall use complex notation $z=x+iy$ for convenience.

In his talk at the Clay Mathematics Institute's Summer School on ``The Global Theory of Minimal Surfaces'' at MSRI in the summer of 2001, Meeks conjectured that there can be at most two such solutions over disjoint domains.

In his paper \cite{S}, Spruck proved this under several stringent
side conditions on the behavior of $u$. Later, without any side conditions, Tkachev \cite{Tk} proved that there could be at most three such disjoint graphs.

We will consider the possible case that there are 3 graphs 
when the domains are simply connected, and provide information
regarding the necessary properties solutions would have in this case.  While these
properties are substantial, they are not sufficient to resolve the conjecture.  

 With $u$ as in (\ref{mineq1}), let $\displaystyle{M(r,u)=\underset{z\in\{|z|=r\}\cap D}{\max} u(z)}$.
 Then
we have the following.

{\bf Theorem 1.}. \emph{Let $u_1,\  u_2,\  u_3$ be solutions to (\ref{mineq1}) over disjoint
simply connected domains $D_1,\ D_2,\ D_3$ respectively.  Then for sufficiently large r,there exist
positive constants $\kappa_1$ and $\kappa_2$ such that }
\begin{equation}
\label{max}
\kappa_1 r<M(r,u_j)<\kappa_2 r\qquad (j=1,2,3).
\end{equation}

We note that for sufficiently large $c,$, by replacing the $u_j$ by $u_j-c$ in the
subdomains $u_j>c$ and using the theorem of Tkachev, we may assume that
each $D_j $ is bounded by a real analytic Jordan arc.

\section{preliminaries}
Let $u$ be a solution to (\ref{mineq1})  over a domain $D$
bounded by an unbounded
differentiable Jordan arc.
We shall  make use of the parametrization of a surface given by $u$ in isothermal coordinates using 
Weierstrass functions $\left( x(\zeta), y (\zeta ), U (\zeta) \right)$
with $\zeta$ in the parameter space.   Our notation will then be
given by 
\begin{equation}
\label{downstairs}
f(\zeta) = x(\zeta) +iy(\zeta).
\end{equation}
Note that we shall use $\zeta =\sigma+i\tau$ in the parameter space, and $z=x+iy$ in $D$.

Then $f(\zeta )$ is univalent and harmonic, and since $D$ is simply connected it can be written in
the form 
\begin{equation}
\label{decomp1}
f(\zeta) = h(\zeta) + \ol{g(\zeta)}
\end{equation}
where $h(\zeta )$ and $g(\zeta)$ are analytic, and
\begin{equation}
\label{Dilatation}
|h'(\zeta)|>|g'(\zeta)|.
\end{equation}
We dismiss the trivial case $g'\equiv 0$ so that
\begin{equation}
\label{decomp2}
U(\zeta )= \pm\  2\Re e \, i\int \sqrt{h'(\zeta )g'(\zeta )}\,d\zeta  .
\end{equation}
(cf. \cite[\S 10.2]{D}).

If the parameter space is taken to be the right half plane $H$, then 
 $z=f(\zeta),$ \newline $ u(f(\zeta))=U(\zeta)$, and $U(\zeta )$ is harmonic and positive in
 $H$ and vanishes on $\partial H$.  Thus,
(cf. \cite[p. 151]{T}), 
\begin{equation}
\label{height}
U(\zeta )= K\,\Re e\, \zeta,
\end{equation}
where $K$ is a positive constant.  This with
(\ref{decomp2}) gives
 
$$
 g'(\zeta) = - \frac{C}{h'(\zeta)}
$$
where $C$ is a positive constant.  The letters $C,\ K,\ B,\ \kappa$ will all be used throughout to
denote constants that may depend on $u$.

By reparametrizing we can make
\begin{equation}
\label{rescale}
U(\zeta ) =2 \Re e\,\zeta \ \ \textrm{and}\ \   g'(\zeta )=-1/h'(\zeta )
\end{equation}
and then the \emph{analytic dilatation} $a(\zeta)$ is
\begin {equation}
\label{dilatation1}
a(\zeta)=\overline{{f_{\overline \zeta}(\zeta)}}/{f_\zeta(\zeta)}=-1/h'(\zeta)^2.
\end{equation}
Furthermore,  from (\ref{Dilatation}) we have, in particular, that
\begin{equation}\label{decomp3}
|h'(\zeta)|=1/|g'(\zeta)|>1.
\end{equation}

For computations in $H$ we shall denote the subsectors by
$$ 
S_\alpha=\{\zeta :|\arg \zeta|<\alpha/2\}.
$$

With these preliminaries we have the following.

{\bf Theorem 2.} \emph{Let $u$ represent one of the solutions in Theorem 1 and parametrized
in $H$ as above.  For any $\alpha <\pi/2$, let $D_\alpha$ denote the image of $S_\alpha$.
Then in $D_\alpha$, the laplacian of $u$ tends to zero as
\begin{equation}
\label{laplacian}
|\Delta u(z)<K/|z| 
\end{equation}
for some $K=K(u,\alpha)$ and $|z|$ sufficiently large }.

\emph{In addition, given $\varepsilon >0$, with $|\alpha|$ sufficiently close to $\pi/2$  we  also have
\begin{equation}
\label{boundary}
u(z)<\varepsilon |z| \qquad z\in\partial D_\alpha
\end{equation}
 for sufficiently large $|z|$.}

The proof of Theorem 2, which will be carried out in \S3-\S10 uses the fact that, with  
(\ref{decomp3}), 
 $\log |h(\zeta)|$ is a positive harmonic function and therefore has a Poisson
representation.

{\bf Lemma A.} \emph{With $h(\zeta)$ as in (\ref{decomp1}) and (\ref{decomp3}),}

\begin{equation}
\label{Poisson1}
\log| h'(\zeta)|=\frac {\sigma} \pi\int_{-\infty}^\infty\frac{\log|h'(is)|ds}{|\zeta -is|^2}  \qquad\zeta\in H.
\end{equation}

{\bf Proof.}  Since by (\ref{decomp3}) $\log|h'|$ is a positive harmonic function in $H$ it can be represented (cf. \cite[p. 149]{T})
$$
\log| h'(\zeta)|=\frac {\sigma} \pi\int_{-\infty}^\infty\frac{\log|h'(is)|ds}{|\zeta -is|^2}+c\,\zeta\qquad c\geq 0  .
$$
Furthermore, by \cite[Theorem IV. 19.]{T}, $\log h'(\zeta)/\zeta\to c$ and 
$h''(\zeta)/h'(\zeta)\to c$ in $S_\epsilon$.  By \cite[Lemma 2]{W1} it follows that $c=0$.

\qed

We shall also need the \emph{conjugate function}, given in $D$ by (e.g.  (\cite[p. 62] {M})
\begin{equation}
\label{conjugate}
v(z)=\int -\frac{u_{w_2}}{\sqrt{1+|\nabla u|}}dw_1+\frac{u_{w_1}}{\sqrt{1+|\nabla u|}}dw_2.
\end{equation}

As in (\ref{rescale}), in the parameter half plane $H$ we may take this as
\begin{equation}
\label{conjugate1}
V(\zeta)=2\Im m\, \zeta.
\end{equation}

\section{Miklykov's method}

  We shall take advantage of the work of V. Mikljukov \cite{M} on the moduli of path
families on minimal surfaces. 

  Let $\ D\ $ be a simply connected
unbounded domain in $\Bbb R^2$.  Let $\ F\ $ denote the surface given by
$\ u(z)$,\quad  $z = x_1 +ix_2 \in D\ $ with 
$$
ds_F^2=(1+u_{x_1}^2)dx_1^2+2u_{x_1}u_{x_2}dx_1dx_2+(1+u_{x_2}^2)dx_2^2
$$ and
$$
d\sigma_F=\sqrt{1+|\nabla u|^2}\,dx_1dx_2
$$ 
the respective length and area
elements for $F$.

For a family $\Gamma\ $ of curves in $\ D\ $ we define the {\it modulus} of
$\Gamma$ in the metric of $\ F\ $ by
$$
\M_F\Gamma =\inf \iint_D \rho^2(z)d\sigma_F,
$$
the $\inf$ being taken over all nonnegative measurable functions $\rho \ $ on 
$\ D\ $ satisfying 
  $$\ds\ \underset{\gamma\in\Gamma}{\inf} \int_\gamma \rho(z)ds_F\ge 1.$$

The utility of the modulus comes from the elementary observation that it is a
conformal invariant (cf. \cite[p.65]{M}).

We shall use estimates on the modulus for path families of curves on a surface
$\ F\ $ given by solutions $\ u(z)\ $ to the minimal surface equation over
domains $\ D\ $ as in (\ref{mineq1}).  Corresponding to this we
also have the conformal parametric description given by (\ref{downstairs}) -(\ref{decomp3}).

Path families $\Gamma$ in $H$ correspond to path families on $F$ which  project to path families
$f(\Gamma )$ in $D$.  By conformal invariance, the modulus may be computed either in $H$, 
or with the surface       
metric in $D$.

Let $D(r,R)=D\cap \{r<|z|<R\}$.
  Fix $a$ and $b$ on $\partial D$. 
 For sufficiently large $\ t>0\ $ let $\ S(t)\ $ be the component of
$\ D\cap\{z:|z|=t\}\ $ separating $\ a\ $ from $\infty\ $ in $\ D$.  Choosing $r$ large
enough so $S(r)$ separates $b$ from $\infty$ in $D$, and $R>r$,  let $T $
 be the subdomain of $\ D\ $ between
$\ S(r)\ $ and $\ S(R)\ $.
Let $\Gamma=\Gamma(r,R)\ $  be a family of
curves in $T$ that join $\ S(r)\ $ and $\ S(R)$.

In what follows, we shall extend the results of \cite{M}  slightly to suit our needs.  We assume
there are sequences $\{r_m\},\ \{r_m'\},\ \{r_m''\}$ tending to $\infty$ such that 
\begin{equation}
\label{pp}
r_m/r_m'\to \infty \textrm{ and } r_m''/r_m\to \infty.
\end{equation}

{\bf Lemma 1.  }  \emph{ With the above notations, if for some $C>0$,
\begin{equation}
\label{initial1}
u(r_mz)/(r_mz)\leq C\quad r_m'\leq |r_mz|\leq r_m''\quad r_mz\in D,
\end{equation}
and
 $M_n=\underset{S(r_n'))}{\max}(|z|/r_n|, r_nu(z/r_n) )$
then
\begin{equation}
\label{4.2a}
 \M_F\Gamma (r_n',r_n'') \leq \frac\pi{2(\log r_n'')^2}\left(\underset{D(r_n',r_n'')}
\iint  \frac{dx_1dx_2}{|z|^2}
)+2\pi^2 +o(1)\right)\ (m\to\infty).
\end{equation}
With the stronger condition
\begin{equation}
\label{initial 1a}
u(r_mz)/(r_mz)\to 0\quad r_m'\leq |r_mz|\leq r_m''\quad r_mz\in D,
\end{equation}
then
\begin{equation}
\label{4.2aa}
\M_F\Gamma (r_n',r_n'') \leq \frac 1{(\log r_n'')^2}\left(\underset{D(r_n',r_n'')}
\iint  \frac{dx_1dx_2}{|z|^2}
)+2\pi^2 +o(1)\right)\ (m\to\infty).
\end{equation}}

{\bf Proof.  }
For notational simplicity we shall replace $z/r_m,\  r_m'/r_m,\ r_m''/r_m,\ u(r_mz)/r_m,$ and $ M_n$ by $z,\ r,\ R,\ u(z),$ and $ M$
respectively.

Following Mikljukov, we choose the density function
$$
\rho(z)=(|z|^2+u^2(z))^{-1/2},
$$
for $z\in B=\overline T\cap \overline{D(r,R)}$ and $\rho(z)=0$ for all the remaining
values $\ z\in D$.  Hence
\begin{equation}
\label{4.3}
\M_F \Gamma \le {\ds\iint_B (|z|^2+u^2(z))^{-1}d\sigma_F\over
\left(\underset{\gamma\in\Gamma}{\inf} \int_\gamma (|z|^2+u^2(z))^{-1/2}ds_F\right)^2}.
\end{equation}

For the denominator in (\ref{4.3}), let $\gamma$ be a curve in $\Gamma$.   Then
\begin{equation}
\label{4.4}
\int_\gamma \rho\, ds_F = \int_\gamma (|z|^2+u^2(z))^{-1/2} ds_F\ge \log\left(1+{R-r\over M}\right).
\end{equation}

To obtain a bound for the numerator we let
\begin{equation}
\label{4.5a}
\eta(z)={1\over |z|}\,\arctan \delta(z),\qquad \delta(z)={u(z)\over |z|},
\end{equation}
so that
\begin{equation}
\label{4.5b}
{1\over |z|^2+u^2(z)}={1\over |z|^2(1+\delta^2(z))}.
\end{equation}

To begin with, we note that
$$
\eta_{x_i}={-x_i\over |z|^3}\,\arctan \delta(z)+{1\over|z|} {1\over 1+\delta^2(z)}\left({u_{x_i}\over |z|}-{x_iu\over |z|^3}\right),
$$
and so
\begin{equation}
\label{4.6}
\sum_{i=1}^2 {u^2_{x_i}\over |z|^2(1+\delta^2(z))\sqrt{1+|\nabla u|^2}}\\ =
 \sum_{i=1}^2 {\eta_{x_i}u_{x_i}\over \sqrt{1+|\nabla u|^2}}
 \end{equation}
 $$
+{\arctan \delta(z)\over |z|^3} \sum_{i=1}^2 {x_iu_{x_i}\over \sqrt{1+|\nabla u|^2}}
+{1\over |z|^3} {\delta(z)\over 1+\delta^2(z)} \sum_{i=1}^2
{x_iu_{x_i}\over \sqrt{1+|\nabla u|^2}}.
$$
Also,
$$
\sqrt{1+|\nabla u|^2}={1\over \sqrt{1+|\nabla u|^2}}+\sum_{i=1}^2
{u_{x_i}^2\over \sqrt{1+|\nabla u|^2}},
$$
and thus 
\begin{equation}
\label{4.7}
\sum_{i=1}^2 {u^2_{x_i}\over |z|^2(1+\delta^2(z))\sqrt{1+|\nabla u|^2}} =
{\sqrt{1+|\nabla u|^2}\over |z|^2(1+\delta^2(z))}
\end{equation}
$$
-{1\over 
|z|^2(1+\delta^2(z))\sqrt{1+|\nabla u|^2}}.
$$

Finally, by Green's Theorem
\begin{equation}
\label{4.8}
\iint_B \sum_{i=1}^2 {\eta_{x_i}u_{x_i}\over\sqrt{1+|\nabla u|^2}}\,dx_1dx_2+
\iint_B \eta \sum_{i=1}^2 {\partial\over\partial x_i} {u_{x_i}\over \sqrt{1+|\nabla u|^2}}\,dx_1dx_2
\end{equation}
$$
= \int_{\partial B} {\eta \nabla u\cdot n\over \sqrt{1+|\nabla u|^2}} |dz|,
$$
with $n$ being the outer unit normal.   
The middle term vanishes because of the minimal surface equation.    

Using (4.5)--(4.7) in (4.8) we then have
%$$
%\sum_{i=1}^2 {x_iu_{x_i}\over |z|\sqrt{1+|\nabla u|^2}}\le
%{|\nabla u|\over \sqrt{1+|\nabla u|^2}}
%$$

$$
\iint_B {d\sigma_F\over|z|^2+u^2(z)} = 
\iint_B {1\over (1+\delta^2(z))\sqrt{1+|\nabla u|^2}}
{dx_1dx_2\over |z|^2}
$$
$$
+ \int_{\partial B} {\eta \nabla u\cdot n\over \sqrt{1+|\nabla u|^2}}\,|dz|
+\iint_B  \arctan \delta(z)\  \sum_{i=1}^2 
{x_iu_{x_i}\over |z|\sqrt{1+|\nabla u|^2}}{dx_1dx_2\over |z|^2}
$$
$$
+ \iint_B {\delta(z)\over 1+\delta^2(z)}\ \sum_{i=1}^2 
{x_iu_{x_i}\over
|z|\sqrt{1+|\nabla u|^2}} {dx_1dx_2\over |z|^2}\\
$$
Since 
$$
\sum_{i=1}^2 {x_iu_{x_i}\over |z|\sqrt{1+|\nabla u|^2}}\le
{|\nabla u|\over \sqrt{1+|\nabla u|^2}},
$$
this gives
$$
\iint_B {d\sigma_F\over|z|^2+u^2(z)} \le 
   \iint_B 
{1\over (1+\delta^2(z))\sqrt{1+|\nabla u|^2}}{dx_1dx_2\over |z|^2}
$$
$$
+\int_{\partial B} {\eta\nabla u\cdot n\over \sqrt{1+|\nabla u|^2}}
{|dz|\over |z|}
+ \iint_B \arctan \delta(z){|\nabla u|\over \sqrt{1+|\nabla u|^2}}
{dx_1dx_2\over |z|^2}
$$
$$
+\iint_B {\delta(z)\over 1+\delta^2(z)}
{|\nabla u|\over \sqrt{1+ |\nabla u|^2}}{dx_1dx_2\over |z|^2}.
$$ 
Since
$$
\ds\int_{\partial B} {\eta\nabla u\cdot n\over \sqrt{1+|\nabla u|^2}}{|dz|\over |z|}\le 
\ds\int_{\partial B} {\eta |\nabla u|\over \sqrt{1+|\nabla u|^2}}{|dz|\over |z|}\le 2\pi^2,
$$
we have
$$
\iint_B {d\sigma_F\over |z|^2+u^2(z)}\le \iint_B \frac1{\sqrt{1+\nabla u|^2}}
 (\frac 1{1+\delta^2(z)}
$$
$$
+|\nabla u|\arctan\delta (z)+
\frac{|\nabla u|\delta (z)}{1+\delta^2(z)})
(1+o(1)) {dx_1dx_2\over |z|^2}
+2\pi^2 \quad (R\to \infty).
$$
By letting $\delta (z)\to \infty$ and using (\ref{initial1}) we obtain
  (\ref {4.2a}).  Using (\ref{initial 1a}) and lettting $\delta\to 0$ we obtain
  (\ref{4.2aa}).
  
   \qed

We now use the conformal invariance of the mapping $H\to F$ as described in
\S3 together with the estimate of Lemma 1.  Continuing the convention of using
 $z,\ r,\ R,\ u(z),$ and $ M$ for $z/r_m,\  r_m'/r_m,\ r_m''/r_m,\ u(r_mz)/r_m,$ and $ M_n$ 
 respectively,  
 with $a,\ b,\ f(\zeta )$ as in \S3,  and continuing  with \cite[p.67]{M}, we take $r>0$ so that
$S(r)$ separates $b$ and $\infty$ in $D$.  For $t\ge r$, let $S^*(t)=f^{-1}(S(t))$ 
so that $\overline{S^*(t)}$ has endpoints on $\partial H$ in the $\zeta$ plane.
 Let  $l(t)$ denote the Jordan curve formed by $\overline{S^*(t)}$ along with its
reflection across  $\partial H$ and $G$ the annular domain between $l(r)$ and $l(R)$.
 Let $\tilde\Gamma(r,R)$ be the family of curves separating $l(r)$ and $l(R)$ in $G$.  Then,
 as in (\cite[p.67]{M}),
since $l(r)$ and $l(R)$ separate 0 ($=f^{-1}(a)$) and 1 ($=f^{-1}(b)$) from $\infty$,
 the modulus (in the Euclidean metric) satisfies \cite [pp. 32,56 and 61 (2.10)]{LV},
$$
\textrm{mod}\,\tilde\Gamma (r,R)\le \frac1{2\pi}{ \log(16(P+1))},
$$
where  
$$
P=\underset{\zeta\in l(R)}{\min} |\zeta |. 
$$
Now let $\Gamma^*(r,R)$ be the curves joining $l(r)$ and $l(R)$ in $G$.  Then $\M\,\Gamma^*(r,R) 
=1/\M\, \tilde\Gamma(r,R)$.    Therefore,
\begin{equation}
\label{4.9}
\textrm{mod}\,\Gamma^* (r,R)\ge {2\pi\over \log(16(m(R)+1))}.
\end{equation}
Let
\begin{equation}
\label{m(t)}
m(t)=\underset{\underset{z\in D}{|z|=t}}{  \min  }|\zeta(z)|.
\end{equation}
Then by (\ref{4.9}), the symmetry principle  
\cite[p.16]{A}, and conformal invariance,
\begin{equation}
\label{4.10}
\textrm{mod}_F\Gamma (r,R)\ge {\pi\over \log(16(m(R)+1))}.
\end{equation}
Thus, (\ref{4.10}) taken together with Lemma 1 yields

{\bf Lemma 2.  } \emph{With the above notation and conventions, then with (\ref{initial1}) we have}
\begin{equation}
%\label{Theorem3a}
\underset{m\to\infty}{\liminf}\frac{\log m(r_m'')}{(\log r _m'')^2}
(\iint_{D(r_m',r_m'')}
(1+o(1)){dx_1dx_2\over |z|^2}+2\pi^2)\ge 2,
\end{equation}

\emph{and if (\ref{initial 1a}) holds, then}
\begin{equation}
%\label{Theorem3aa}
\underset{m\to\infty}{\liminf}\frac{\log m(r_m'')}{(\log r _m'')^2}
(\iint_{D(r_m',r_m'')}
(1+o(1)){dx_1dx_2\over |z|^2}+2\pi^2)\ge \pi.
\end{equation}

In using Lemma 2, we shall need an estimate on the quantity $m(t)$.
\vskip .1truein
{\bf Lemma 3.}  \emph{Let $u(z)$ be as in Lemma 1 and Lemma 2.   
then $M(t)\geq m(t)$. }  
\vskip .1truein
{\bf Proof.  }If $\zeta_0$ gives the minimum in (\ref{m(t)}), then $m(t)=|\zeta_0|$.  If $z_1$ is a point
on $|z|=t$ where $M(t)$ is achieved and $f(\zeta_1)=z_1$, then by (\ref{rescale}), $u(z_1)
=u(f(\zeta_1))=\re \,\zeta_1\geq
|\zeta_0|$ or, in other words, $M(t)\geq m(t)$.   \qed

\section{growth lemmas}

For a given $u(\zeta)$ as in the statement of Theorem 1 over a domain $D$, we define $\Theta(r)
=\Theta(u,r)$ to be the angular
measure of the set $D\cap\{|\zeta=r\}$.

{\bf Lemma 4.} \emph{Let $D$ be an unbounded simply connected 
domain with with $\Theta(r_n)\leq\beta<\pi$
    for an increasing sequence $r_n\to\infty$.
    If $u$ is a solution to the boundary value problem (\ref{mineq1}). 
then there exists a sequence $r_n'<r_n$ such that $r_{n}/r_n'\leq B$ for some $B>0$. 
and a constant $K=K(\beta,B )$, such that }
$$
u(z)\leq Kr_n\qquad z\in D\cap\{|z|=r_n'\}.
$$

\emph{In case $r_{n+1}/r_{n}$ remains bounded, it then follows readily that
\begin {equation}
\label{ugrowth}
 M(r,u)<Cr
 \end{equation}
 for some $C>0$ and all large $r$.}

{\bf Proof.  } Let $\gamma$ be an arc on the unit disk having length 
$\beta$ 
and $\ell$ be the line segment joining the endpoints of $\beta$. Let $V$ be the subset of the disk having its boundary as $\ell$ and the complement of $\beta$ with respect
to the boundary of the disk.

Extend $\ell$ symmetrically to length 1, and let $\tilde\ell$ be the line formed by translating
the extended $\ell$ parallel to itself 10 units in the direction of the origin. Now construct
a rectangle $R$ by joining the endpoints of the extended $\ell$ to the endpoints of 
$\tilde\ell$.  Then $V\subset R$ and it follows from \cite[p. 322]{JS}
that there exists a solution $v(z)$ to the minimal surface equation in $R$ having boundary values $+\infty$ on the extended $\ell$ and $\tilde\ell$, and 0 on the other two sides.

 Let $\eta$ be the distance from the origin to $\ell$
and $\Delta=\{|z|<\eta/2\}$.  Let $M$ denote the maximum of $v(z)$ in $\Delta$.

Now, for any given $r_n$, the function $u(r_nz)/r_n$ satisfies the minimal surface equation
and we consider any of the components of the set where $u(r_nz)>0$ in
 $|z|<1$.  In any such component $u(r_nz)/r_n$
is bounded by the corresponding $v(z)$ with the arc $\gamma$ referred to above, being
an arc containing the portion of the component on the boundary.  Thus, if the component
extends to $\Delta$, then for those points in $\Delta$ we have
$u(r_nz)/r_n\leq M$.  Applying this to all such components we then get $u(z)
\leq Mr_n$ for $z\in D\cap\{|z|=\eta\, r_n/2\}$ from which the result follows.

\qed  

{\bf Lemma 5.}\emph { Let $u$ satisfy (\ref{mineq1}) and
assume that for sequences
$ \{r_n'\}$ and  $\{r_n''\}$ satisfying $r_n''/r_n'\to\infty$, and that  (\ref{ugrowth}) holds for
$r_n'<|z|<r_n''$.  Then,}
\begin{equation}
\label{i}
\underset{n\to\infty}{\liminf}\frac 1 {\log r_n''}(1+o(1))\iint_{D(r_n',r_n'')}\frac{drd\theta}{r}\geq 2,
\end{equation}

\emph{and if (\ref{initial 1a}) holds, then}
\begin{equation}
\label{ii}
\underset{n\to\infty}{\liminf}\frac 1 {\log r_n''}(1+o(1))\iint_{D(r_n',r_n'')}\frac{drd\theta}{r}\geq \pi,
\end{equation}

{\bf Proof.} This follows directly from Lemma 3 and Lemma 4.\qed
%With $\{r_n'\}$ and $\{r_n''\}$
%as above,
 %and $f$ as in \S2, 
%by (\ref{rescale}) it follows that $U(m(r_n''))=2m(r_n'')$.  By the maximum principle, 
%$|f(\zeta))|\leq r_n''$ if $|\zeta|\leq m(r_n'')$.  Then, using (\ref{ugrowth}),  we have 
f%or some $\zeta$ with $|\zeta|=m(r_n'')$,
%$$
%m(r_n'')=U(m(r_n'')/2=U(f(\zeta))/2\leq Kr_n''
%$$
%for a $K>0$ from which (\ref{i}) follows from (\ref{Theorem3a}).\qed

\section{Proof of (\ref{max})\phantom{xxxxxxxxxxx}}
%{\bf Lemma 4. } \emph{Let $u_1,\ u_2, u_3$ be solutions to the minimal surface equation as in Theorem 1, with corresponding
%domains $D_1,\ D_2,\ D_3$.  Then each $u_j$ satisfies (\ref{ugrowth}) in its domain $D_j$.}
 We first prove the right side of (\ref{max}). 
 
 Suppose first that there is a sequence $\{r_n'\}$ tend to $\infty$ such that 
$$
M(u_1,r_n')/r_n'=a_n\to\infty\quad\textrm{and}\quad M(u_2,r_n')/r_n'=b_n\to\infty,
$$
and take $d_n=\min (a_n,b_n)$.  Then, if $r_n''=r_n'/d_n^{3/4}$, we have 
\begin{equation}
\label{infinite}
M(u_1,r)/r\to\infty\quad\textrm{and}\quad M(u_2,r)/r\to\infty
\end{equation}
for $r_n'<r<r_n''$ as $n\to\infty$.  It follows from Lemma 4 that for $\varepsilon>0$, $\Theta (u_1,r)
>\pi-\varepsilon$ and $\Theta(u_2,r)>\pi-\varepsilon$ for $r\in I_n=(r_n'^{1/4},r_n'^{1/2})$ and $n$ large.   This 
means that $\Theta(u_3,r)<\pi/2$ on $I_n$ so that repeated use of Lemma 4 applies, and hence so does 
Lemma 5.  However, this is inconsistent with the fact that $\Theta(u_1,r)+\Theta(u_2,r)>2\pi-2\varepsilon$
for this interval.

This limits the possibility down to only one of the solutions, say $u_1$ having  infinite sequence
$\{r_n'\}$ and $\{r_n''\}$ as above tending to $\infty$ such that $u_1(z)/|z|\to\infty$ for $r_n'<|z|<r_n''$, but
$u_2$ and $u_3$ satisfy (\ref{ugrowth}) on these intervals.
 This  possibility then gives $\Theta(u_1,r)  
>\pi-\varepsilon$ on $I_n$
and $u_2$ and $u_3$ satisfying (\ref{i}) which again is impossible.
 Thus,  each $u_j$ satisfies (\ref{ugrowth}) in $D_j$.

Having established the right side of (\ref{max}), we focus on one of the solutions,
say $u_1$, and use the fact that $u_2$ and $u_3$ satisfy (\ref{ugrowth}) and
hence also (\ref{i}).

Now if there exists a sequence $r_n\to\infty$  on which $M(u_1,r_n)/r_n=\delta_n
\to 0$, then $M(u_1,r)/r\to 0$ on the interval $[\sqrt{\delta_n}r_n,r_n]$ so using
this interval for $u_1$ in (\ref{ii}) and also for $u_2$ and $u_3$ in (\ref{i}) we
obtain a contradiction thus proving the left side of (\ref{max}).\qed

\section{the three cases}
In what follows we assume 3 disjoint solutions $u_1,\ u_2, \ u_3$ as in Theorem 1.
By \S 5 we may
apply Lemma 5 to $u_2$ and $u_3$.  As a consequence, it must then be
 that for the sequences $\{r_n'\}$ and $\{r_n''\}$ with $r_n''/r_n'\to\infty$,  $u_1$ in
 its domain $D_1$  must satisfy
 \begin{equation}
 \label{gamma1}
 \underset{n\to\infty}{\limsup}\frac 1 {\log r_n''}(1+o(1))\iint_{D_1(r_n',r_n'')}\frac{drd\theta}{r}\leq
\gamma
\end{equation}
where
\begin{equation}
\label{gamma2}
\gamma=2\pi -4.
\end{equation}

Since we will focus on $u_1$,
  for convenience we will suppress the subscript so $u(z)=u_1(z)$ and $D=D_1$.

 For notational convenience we may rescale the 3 solutions $u_j(cz)/c\ j= 1,2,3$
with $c>0$ chosen so that the $f(\zeta)$ corresponding to $u$ as in \S 2 satisfies 
(\ref{decomp3}).

%The strategy  will be to analyze   (\ref{Poisson1}) corresponding to this $u$ in  sectors,
%where we now use the notation 
%$$
%S_\epsilon =\{-\pi/2 +\varepsilon <\arg \zeta<\pi/2 -\varepsilon\},
%$$
 %and $\varepsilon>0$ is fixed arbitrarily small.   
 We will consider 3 possibilities: 
\bigskip

Case 1.  $\displaystyle{\underset{t\to\infty}{\limsup}\frac 1 t\int_{-t}^t\log|h'(is)|ds=\infty}$

Case 2.  $\displaystyle{\underset{t\to\infty}{\liminf}\frac 1 t\int_{-t}^t\log|h'(is)|ds=0}$

Case 3.  $0<\displaystyle{\underset{t\to\infty}{\liminf}\frac 1 t\int_{-t}^t\log|h'(is)|ds\leq
\underset{t\to\infty}{\limsup}\frac 1 t\int_{-t}^t\log|h'(is)|ds<\infty}$
\eject
Henceforth we shall modify our notation, taking
$$
S_\epsilon =\{-\pi/2 +\varepsilon <\arg \zeta<\pi/2 -\varepsilon\},
$$
 and $\varepsilon>0$ is fixed arbitrarily small. 
\section{Quasiconformal mappings}
We shall have occasion to view the harmonic mapping described in \S 2 as a  quasiconformal mapping.  A one to one sense preserving mapping  $f$ in a domain $D$ is quasiconformal, if its \emph {dilatation} $\delta = \delta (\zeta )$ defined by  (cf.(\cite[p. 5]{A}))
\begin{equation}
\label{dilatation}
\delta =\left |\frac{f_{\overline{\zeta}}}{f_\zeta}\right|
\end{equation}
satisfies
\begin{equation}
\label{quasiconformal}
\underset{\zeta \in D}{\sup}\,\delta(\zeta)<1.
\end{equation}
The dilatation is a conformal invariant, and the inverse mapping has the same dilatation at corresponding
points \cite[p. 9]{A}.

We shall need a modification of the Ahlfors distortion theorem.
In the classical setting we have a simply connected region $G$ with accessible boundary
points $Z_1=X_1+iY_1$ and $Z_2=X_2+iY_2$ 
$Z_1$ and $Z_2$.   We assume that $-\infty\leq X_1=\inf \Re e\,z$ for $z\in G$ 
and $\infty\geq X_2=\sup \Re e\, z$ for $z\in G$. We consider 
 $z=x+iy$ in $G$ with crosscuts 
 $Q_x$  separating $Z_1$ and $Z_2$  in $G$ (See (\cite[pp. 94-95]{Ne}) for more details).  Let $\Theta
(x)$ be the length of $Q_x$. Let $w(z)=u(z)+iv(z)$ be a conformal mapping of $G$ onto the strip
$\{|v|<a/2\}$  such that $Z_1$ corresponds to $-\infty$ and $Z_2$ to $+\infty$.  If $u_1(x)$ denotes the 
smallest value on the crosscut and $u_2(x)$ the largest, then
the classical distortion theorem is as follows.

{\bf Theorem  A. } \emph{If
$$
\int_{x_1}^{x_2} \frac{dx}{\Theta(x)}>2,
$$
then 
$$
u_1(x_2)-u_2(x_1)\geq a\int_{x_1}^{x_2} \frac 1 {\Theta(x)}-4a.
$$
}
For our purposes, the strip $\Sigma_\varepsilon$ will be the (principal branch) logarithmic image of $S_\varepsilon$ in the $w=\mu+i\nu$ plane and $G$ will be  
 the image of a fixed branch of $\log f(S_\varepsilon)$ in the $z=x+iy$ plane  
with $f$ as in \S2. 

Let $w(z)=\log(f^{-1}(e^z))$ for the principal branch of log which then has the same dilatation as $f$
  at corresponding points.  Then \cite[Lemma 1]{W3},

{\bf Lemma B.}  \emph{With the above notations,  
let $R$
be a rectangle in the $\mu +i\nu$ plane
$$
R=\Sigma_\varepsilon\cap\ \{p \leq \mu\leq q\}\quad 0<p<q,
$$
and suppose that  $w(z)$ has dilatation less than $\delta_0$ in $w^{-1}(R)$. 
Then for $x_1+iy_1$ and $x_2+iy_2$ in $w^{-1}(R)$ ($x_1<x_2$) and $a=\pi-\varepsilon$, if
\begin{equation}
\label{Ahlfors1}
 \int_{x_1}^{x_2} \frac{dx}{\Theta(x)}>2\frac{1+\delta_0}{1-\delta_0}, 
\end{equation}
}
then
\begin{equation}
\label{Ahlfors4}
\mu_1(x_2)-\mu_2(x_1)\geq a\frac{(1-\delta_0)}{(1+\delta_0)}\int_{x_1}^{x_2}\frac{dx}{\Theta(x)} 
-4a.
\end{equation}

\section{Case 1}

 Let $I_n = (2^{n-1},2^{n})$ for $n= 1, 2.....$ 
 and $I_n=(-2^{|n|}, -2^{|n|-1})$
  for $n=-1,-2,....$, and 
 $I_0=(0,1)$.  
 Suppose that there exists a subsequence $I_{n_k}$ and $\alpha_{n_k}\to\infty$ such that 
 $$
 \int_{I_{n_k}}\log|h'(is)|ds=\alpha_{n_k}2^{n_k}.
$$
By (\ref{Poisson1}),  for $\zeta\in H$,
\begin{equation}
\label{logh'}
\log|h'(\zeta)|=\frac \sigma \pi\int_{-\infty}^\infty\frac{\log|h'(is)|ds}{|\zeta -is|^2}.
\end{equation}
Using positive $n_k$ for notational convenience, from this it follows that as $n_k\to\infty$,
$$
\log|h'(\zeta)|\geq  C \frac {2^{n_k}}{\sqrt{\alpha_{n_k}}}\frac{\alpha^{n_k}2^{n_k}}{2^{2n_k}}\to\infty\quad
  \zeta\in S_\epsilon\cap
  \{2^{n_k}/\sqrt{\alpha_{n_k}}
  \leq|\zeta |\leq 2^{n_k}\}.
$$
for some constant $C$.

  Summarizing this, we obtain sequences $\{r_n'\}$ and $\{r_n''\}$ both tending to $\infty$ such that 
  $r_n''/r_n'\to\infty$, and with $T_n=S_\epsilon\cap\{r_n'\leq|\zeta|\leq r_n''\}$
  \begin{equation}
  \label{h'large}
  |h'(\zeta )|\to\infty\quad  \zeta\in \overline T_n\ \ n\to\infty.
  \end{equation}
  In order to apply (\ref{h'large}), we consider $f$ as a quasiconformal mapping in $T_n$. As such,
  by (\ref{decomp1}), (\ref{dilatation1}), and (\ref{dilatation}) this becomes $\delta=1/|h'|^2 
  \to 0$
 in $T_n$  by (\ref{h'large}).
  
    Let $w(z)=\log(f^{-1}(e^z))$ which has the same dilatation as $f$
  at corresponding points.  Let $G$ denote the logarithmic image of $f(S_\epsilon )$.  Then $w$ maps $G$ into
  the strip \newline
    $(\pi - \varepsilon)/2 <\Im m \,w<(\pi + \varepsilon )/2$.  Let $R_n$ be the rectangle corresponding
  to $T_n$ in the logarithmic plane with $s_{1,n}$ denoting the left side and $s_{2,n}$ the right side of $R_n$.
  
  Next we take the subsets $\{G_n\}$ of $G$ corresponding  to the preimages $w^{-1}(R_n)$, and take
  $$
  x_{1,n}=\underset{z\in w^{-1}(s_{1,n})}{\sup}\Re e\,z\ \ \textrm{and}\ \ 
    x_{2,n}=\underset{z\in w^{-1}(s_{2,n})}{\inf}\Re e\,z.
  $$
We note that 
\begin{equation}
\label{mu}
\begin{aligned}
\log|f(\zeta)|\leq x_{1,n}\ \ 
 (\zeta\in S_\varepsilon
\cap\{|\zeta|=r_n'\}\\ 
\log|f(\zeta)|\geq x_{2,n}\ \ 
 (\zeta\in S_\varepsilon
\cap\{|\zeta|=r_n''\}
\end{aligned}
\end{equation}

 Now,

 $$
 \int_{\log r_n'}^{\log r_n''}\frac {dx}{\Theta(x)}>\frac 1{2\pi}\log\frac{r_n''}{r_n'}
 $$
 which assures the condition (\ref{Ahlfors1}).
 
 Since
 $$
 \int_{\log r_n'}^{\log r_n''}\frac{dx}{\Theta(x)}
  \int_{\log r_n'}^{\log r_n''}{\Theta (x) dx}\geq(\log\frac{r_n''}{r_n'})^2,
 $$
from (\ref{gamma1})we have 
$$
\int_{\log r_n'}^{\log r_n''}\frac {dx}{\Theta(x)}>\frac {\log{(r_n''}(1+o(1))}{\gamma}
\qquad(n\to\infty).
$$
Using Lemma B we then have
$$
x_{2,n}-x_{1,n}>\frac{\pi-\varepsilon}{\gamma}\log (r_n''(1+o(1)))\qquad(n\to\infty).
$$
Combining (\ref{mu}) with Theorem 1 and (\ref{gamma2}) we then obtain a
contradiction which shows that Case 1 is not possible.

 \section {Case 2}
 
  Having eliminated Case 1, it follows that 
\begin{equation}
\label{Case1}
\int_{-s}^{s}\log|h'(it)|dt<Cs
\end{equation}
for some constant $C>0$.

For Case 2 we may assume that there is a sequence $\{s_n\}$  tending to $\infty$ such that 

 \begin{equation}
\label{beta1}
 \frac{\int_{-s_n}^{s_n} \log|h'(it)|dt}{s_n}\leq\beta_n\to 0,\quad \beta_ns_n\to\infty.
 \end{equation}

% We define
% corresponding intervals $\{I_n\}$ and $ \{J_n\}$ given by
% $$
%I_n=[\beta_n^{3/4}s_n,s_n]
%%and
% $$
%J_n=[\beta_n^{1/2}s_n,\beta_n^{1/4}s_n].
%$$
 
 By (\ref{Case1}) an (\ref{Poisson1}) we have
 $$
 \log|h'(\sigma)|=
 \frac{ \sigma}{ \pi}\int_{-\infty}^{\infty}\frac{\log|h'(is)|ds}{\sigma^2+s^2}
\leq \sigma  \frac{Cs}{\sigma^2+s^2}
 \bigg |_{0}^{\infty}+C\sigma\int_{0}^{\infty}
 \frac {s^2 ds}{(\sigma^2+s^2)^2}.
$$
From this it follows that
$$
\log|h'(\sigma)|<C
$$
 and consequently
\begin{equation}
\label{prelim}
|h(\sigma|<e^C\sigma
\end{equation}
for some constant $C$.

We now break the integral in (\ref{Poisson1}) into parts and focus on the values of $\log |h'|$ on intervals of the
ray $\zeta=\sigma$ with 
$$
\sigma\in J_n=[\beta_n^{1/2}s_n,\beta_n^{1/4}s_n].
$$

Now,
$$
\frac \sigma \pi\int_{-\beta^{3/4}s_n}^{\beta_n^{3/4}s_n}\frac{\log|h'(is)|ds}{\sigma^2+s^2}
\leq C\sigma  \frac{s}{\sigma^2+s^2}
 \bigg |_{0}^{\beta_n^{3/4}s_n}+C\sigma\int_{0}^{\beta_n^{3/4}s_n}
 \frac {s^2 ds}{(\sigma^2+s^2)^2}.
$$

%In the interval $J_n$, the first term on the right is bounded by
%$$
%\frac{C\beta_nt^2_n}{\beta_n^{1/2}t_n^2+\beta_n^{3/2}t_n^2}
%$$
%which tends to $0$.
Simple computations with $\sigma$ on the interval $J_n$ show that the terms on the
right tend to $0$. 

The same approach on $(-\infty,-s_n)\cup(s_n,\infty)$ gives 
$$
\frac \sigma \pi\int_{(-\infty,-s_n)\cup(s_n,\infty)}\frac{\log|h'(is)|ds}{(\sigma^2+s^2)^2}\to 0
$$
for $\sigma\in J_n $ as $n\to\infty$.

Finally, using (\ref{beta1}) we have
 $$
\frac \sigma \pi\int_{(-s_n,-\beta_n^{3/4}s_n)\cup(\beta_n^{3/4}s_n,s_n)}\frac{\log|h'(is)|ds}{\sigma^2+s^2}
\leq 
C \sigma  \frac{\beta_ns_n}{\sigma^2+s^2}
 \bigg |_{\beta_n^{3/4}s_n}^{s_n}+C\beta_ns_n\sigma\int_{\beta_n^{3/4}s_n}^{s_n}
 \frac {s ds}{(\sigma^2+s^2)^2}
 $$
 
where  terms on the right again tend to $0$ for $\sigma\in J_n$.

 Thus, we have
\begin{equation}
\label{realpart}
\log|h'(\sigma)|\to 0\qquad \sigma\to\infty,\ \sigma\in J_n.
\end{equation}

Let $\log h'(\sigma )=\varphi(\sigma)+i\psi(\sigma)$ with $\varphi$ and $\psi$ real, and
by (\ref{realpart}) $\varphi(\sigma)\to 0$ as $\sigma\to\infty$ in $J_n$. 
Then from (\ref{decomp1}) and  (\ref{rescale}), on the intervals $J_n$
$$
|f(\beta_n^{1/4}s_n)-f(\beta_n^{1/2}s_n)|=|\int_{\beta_n^{1/2}s_n}^{\beta_n^{1/4}s_n}
(h'(t)-1/\bar{h'}(t))d\sigma|
$$
$$
=|\int_{\beta_n^{1/2}s_n}^{\beta_n^{1/4}s_n} e^{\varphi(\sigma)+i\psi(\sigma)}
-e^{-\varphi(\sigma)+i\psi(\sigma)}d\sigma|
$$
$$
\leq \int_{\beta_n^{1/2}s_n}^{\beta_n^{1/4}s_n}| e^{\varphi(\sigma)+i\psi(\sigma)}
-e^{-\varphi(\sigma)+i\psi(\sigma)}|d\sigma
\leq  \int_{\beta_n^{1/2}s_n}^{\beta_n^{1/4}s_n}\varepsilon_nd\sigma
$$
where  $\varepsilon_n\to 0$ as $n\to \infty$. 
Using this and (\ref{prelim}) we deduce that
$$ 
u(\beta_n^{1/4}s_n)/|f(\beta_n^{1/4}s_n)|>\beta_n^{1/4}/(C\beta_n^{1/2}+\varepsilon_n\beta_n^{1/4})
$$
which contradicts Theorem 1.  Thus, Case 2 is complete.

\section{Case 3}

Turning now to Theorem 2, we first have

{\bf Lemma 6.} \emph{Let $f$ and $U$ be as in \S 2, and assume the conditions of Case 3 are satisfied.
Then, given $0<\varepsilon < \pi/2$, there exists $C_2>0$ (depending on $U$, but independent of
$\varepsilon$) such that}
$$
|h'(\zeta)|\geq C_2(1+o(1))\exp (\varepsilon)\qquad (\zeta\to\infty,\ \zeta\in S_\varepsilon).
$$
%{\bf Lemma 5.} \emph{Let $f$ and $U$ be as in \S 2, and assume the conditions of Case 3 are satisfied.
%Then,  there exists $C>0$ such that if $\zeta=|\zeta| e^{i\alpha}\in H$,}
%$$
%|h'(\zeta)|\geq C(1+o(1))\exp (\tan(\frac \pi 2-\alpha)\qquad (\zeta\to\infty,\ \zeta\in H).
%$$
{\bf Proof.} Using (\ref{Poisson1}) and moving the ``mass'' to the positive $y$ axis
 we may write 
\begin{equation}
\label{Poisson2}
\log| h'(\zeta)|\geq \frac {\sigma} \pi\int_{0}^\infty\frac{\Psi (s)}
{\sigma^2 + (s+t\sigma)^2}ds \qquad(\zeta\in S_\varepsilon),
\end{equation}
where $\Psi(s)=\log|h'(is)+\log|h'(-is)|$, $|\zeta|=|\sigma + it\sigma|$ and 
$\displaystyle{t=\tan(\frac\pi 2 - \varepsilon)}$.

Integrating (\ref{Poisson2}) by parts, the right side of (\ref{Poisson2}) becomes
\begin{equation}
\label{Poisson3}\
\frac \sigma\pi\left(\frac{\int_0^s\Psi(x)dx}{\sigma^2+(s+t\sigma)^2}\bigg |_0^\infty+2\int_0^\infty
\frac{(s+t\sigma)\int_0^s\Psi(x)dx}{(\sigma^2+(s+t\sigma)^2)^2}ds\right).
\end{equation}
From the elimination of Case 1, it follows that the first term in (\ref{Poisson3}) vanishes.

Regarding the second term, using the elimination of Case 2 we can write
$$
\frac{2\sigma}\pi\int_0^\infty
\frac{(s+t\sigma)\int_0^s\Psi(x)dx}{(\sigma^2+(s+\sigma t)^2)^2}ds \geq (1+o(1))\frac{2\sigma}\pi
\int_0^\infty\frac{(s+t\sigma)Cs}{(\sigma^2+(s+t\sigma)^2)^2}ds
$$
as $\sigma\to\infty$.

Since
$$
\sigma\int_0^\infty\frac{(s+t\sigma)s}{(\sigma^2+(s+t\sigma)^2)^2}ds=
\frac 1 2\left (\arctan(s/\sigma+t)-\frac{s\sigma}{\sigma^2+(s+t\sigma)^2}\right )\bigg |_{s=0}^\infty
=\pi/2-\arctan ( t)
$$
Lemma 6 follows.
\qed
\eject

In order to estimate $\Delta u$ we use estimates from
\cite{W1}.

{\bf Theorem E.  } \emph{There exists a constant $K$ depending on $u$ such that, if $u$ is as in
(\ref{mineq1}), and $C>0$, then the curvature $\kappa=\kappa (C)$ 
of the level set $u=C$ satisfies  the 
inequaltiy} 
\begin{equation}
\label{curvature1} |\kappa|\leq \frac KC.
\end{equation}

\vskip .1truein
{\bf Proof of Theorem 2.}. 
  The  curvature of a 
level curve through a point is given by (cf. \cite[p72]{Gray})
\begin{equation}
\label{curvature2}
\kappa = \frac {u_y^2u_{xx}+u_x^2u_{yy}-2u_xu_yu_{xy}}{|\nabla u|^3}
\end{equation}
and by (\ref{mineq1}),
\begin{equation}
\label{Delta}
\Delta u=-u_y^2u_{xx}-u_x^2u_{yy}+2u_xu_yu_{xy}.
\end{equation}

Furthermore in $S_\varepsilon$ we have by \cite[(p.105)]{O} and (\ref{dilatation1}),
$$
|1/h'(\zeta)|^2= |a(\zeta )| =\frac{\sqrt{1+|\nabla u(z)|^2}-1}{\sqrt{1+|\nabla u(z)|^2}+1} ,
$$
and hence,
\begin{equation}
\label{grad1}
|\nabla u|^2=\left(\frac{|h'|^2+1}{|h'|^2-1}\right)^2-1 <\frac 4{(|h'|-1)^2}
\end{equation}

The statement (\ref{laplacian}) then follows from (\ref{curvature2}), (\ref{Delta}), (\ref{grad1}) and Lemma 6.

Regarding (\ref{boundary}),
 starting with (\ref{max}), we observe that for large $z\in f(S_{\pi/4})$ and $v(z)$ as in
 (\ref{conjugate}), then
$v(z)<u(z)<\kappa_2|z|$. For $\alpha>\pi/4$, $v(z)>u(z)$ in $f(S_{\pi/4})$.  Since 
$|\nabla v(z)|<1$ we have  $|v(z)|<\kappa_2|z|+2\pi |z|$.  Thus, when $\tan\,\alpha
=(\kappa_2+2\pi)/ {\varepsilon}$, (\ref{boundary}) holds.

\qed
\eject

\bibliographystyle{amsplain}

\end{document}